\documentstyle{article}

\newtheorem{theorem}{Theorem}[section]
\newtheorem{lemma}[theorem]{Lemma}
\newtheorem{proposition}[theorem]{Proposition}
\newtheorem{corollary}[theorem]{Corollary}
\newtheorem{definition}[theorem]{Definition}

\newtheorem{example}[theorem]{Example}
\newcommand{\GK}{\mbox{GK}}
\newcommand{\ch}{\mbox{char}}
\newcommand{\ann}{\mbox{ann}}
\newcommand{\supp}{\mbox{supp}}
\newcommand{\lin}{\mbox{lin}}

\begin{document}

\title{Quadratic algebras of skew type and the underlying
        semigroups\thanks{Work supported in part
        by  Onderzoeksraad of Vrije Universiteit Brussel,
        Fonds voor Wetenschappelijk Onderzoek (Belgium),
        and KBN research grant
      2P03A 030 18 (Poland).\newline MSC Codes: 16P40, 20M25,  16S15 and 16S36} }
\author{T. Gateva-Ivanova, Eric Jespers and Jan Okni\'{n}ski}
\date{}
\maketitle

\begin{abstract}
We consider algebras over a field $K$ defined by a presentation
\newline $K \langle x_{1},\ldots ,x_{n} : R \rangle $, where $R$
consists of $n\choose 2$ square-free relations of the form
$x_{i}x_{j}=x_{k}x_{l}$ with every monomial $x_{i}x_{j}, i\neq j$,
appearing in one of the relations. Certain sufficient conditions
for the algebra to be noetherian and PI are determined. For this,
we prove more generally that right noetherian algebras of finite
Gelfand-Kirillov dimension defined by homogeneous relations
satisfy a polynomial identity. The structure of the underlying
monoid, defined by the same presentation, is described. This is
used to derive information on the prime radical and minimal prime
ideals. Some examples are described in detail. Earlier, Etingof,
Schedler and Soloviev, Gateva-Ivanova and Van den Bergh, and the
authors considered special classes of such algebras in the
contexts of noetherian algebras, Gr\"{o}bner bases, finitely
generated solvable groups, semigroup algebras, and set theoretic
solutions of the Yang-Baxter equation.
\end{abstract}

\section{Introduction} \label{Section1}

We consider finitely generated monoids with a monoid presentation
of the form $$S=\langle x_{1},x_{2},\ldots ,x_{n} \, | \,
x_{i}x_{j}=x_{k}x_{l} \rangle $$ with ${n\choose 2}$ relations,
where $i\neq j, k\neq l$ and every product $x_{p}x_{q}$ with
$p\neq q$ appears in one of the relations. So each $x_{p}x_{q}$
appears in exactly one relation. We call such an $S$ a semigroup
of skew type. Special classes of monoids of this type, and
algebras defined by the same presentations, arise in a natural way
from the study of set-theoretic solutions of the Yang-Baxter
equation and independently from certain problems in the theory of
associative algebras, \cite{eting},\cite{gateva},\cite{binom}.
These algebras turn out to have very nice properties. In
particular, they have finite global dimension, satisfy the
Auslander regularity condition and they are Cohen-Macaulay
\cite{gat-van}. Reasons and tools for dealing with these
properties came from the study of homological properties of
Sklyanin algebras by Tate and Van den Bergh \cite{tate}.

The above mentioned special classes of semigroups also satisfy the
requirement that $i>j$, $k<l$, $i>k$, $j<l$ for each of the
relations $x_{i}x_{j}=x_{k}x_{l}$, and surprisingly, they define
submonoids of torsion-free abelian-by-finite groups. In
particular, for any field $K$, the semigroup algebra $K[S]$ is a
domain that satisfies a polynomial identity and one also shows
that it is left and right noetherian. Notice that, under this
additional assumption on the relations, every element of $S$ can
be written in the form $x_{1}^{k_{1}}\cdots x_{n}^{k_{n}}$ for
some non-negative integers $k_{i}$. In particular, the
Gelfand-Kirillov dimension of $K[S]$, denoted by $\GK (K[S])$,
does not exceed $n$. We note that some other (but related) types
of algebras defined by quadratic relations have been investigated,
see for example \cite{gateva},\cite{gat-van},\cite{kramer}.

Our aim is to study the noetherian property of algebras $K[S]$ of
skew type, its relation to the growth and the PI-property, and the
role of the minimal prime ideals with respect to the least
cancellative congruence on $S$. This is motivated by the results
on algebras of binomial semigroups, where the height one primes
turned out to be crucial for the properties of the algebra,
\cite{binom}.

Our main result asserts that $K[S]$ is a noetherian PI algebra for
a wide class of semigroups of skew type. A combinatorial approach
allows us to derive a rich structural information on $S$. This is
of independent interest and becomes the main tool in the proof. As
an intermediate step we prove the following general result.
Suppose $A$ is a unitary $K$-algebra defined via a presentation
$K\langle x_{1}, \cdots , x_{n} : R \rangle$, where $R$ consists
of relations of the type $u=v$ with $u$ and $v$ words of equal
length in the generators. If $A$ is right noetherian and of finite
Gelfand-Kirillov dimension then $A$ satisfies a polynomial
identity.

\section{Cyclic condition}

We start with a combinatorial condition that allows us to build
several examples of noetherian PI-algebras $K[S]$. If $S$ is a
monoid and $Z\subseteq S$ then we denote by $\langle Z \rangle$
the submonoid generated by $Z$.

We say that a monoid $S$ generated by a finite set $X$ satisfies
the cyclic condition (C) if for every pair $ x,y \in X$ there
exist elements $x=x_1, x_2, \ldots ,x_k$, $y'\in X$ such that
\begin{eqnarray*} \label{cyc}
 yx&=&x_2y'\\ yx_2&=&x_3y'\\ &\vdots & \\ yx_k&=&xy'
 \end{eqnarray*}
It is shown in \cite{gat-van} (see also \cite{binom}) that
binomial semigroups satisfy the cyclic condition. We show that the
cyclic condition is symmetric.

\begin{proposition}       \label{full-cyclic}
Let $S=<X;R>$ be a semigroup of skew type. Assume $S$ satisfies
the cyclic condition. Then the full cyclic condition (FC) holds in
$S$, that is, for any pair $x,y\in X$, there exist two sequences:
$x=x_1, x_2, \ldots , x_k$ and $y = y_1, y_2,\ldots ,y_p $ in $X$
such that $$y_1x_1 =x_2y_2 , y_1x_2=x_3y_2,\ldots , y_1x_k=x_1y_2,
$$ $$y_2x_1 =x_2y_3, y_2x_2=x_3y_3,\ldots , y_2x_k=x_1y_3, $$
$$\vdots $$ $$y_px_1 =x_2y_1, y_px_2=x_3y_1,\ldots ,
y_px_k=x_1y_1.$$
\end{proposition}

We call this a cycle  of type $k\times p$.

\begin{lemma} Under the hypothesis of Proposition~\ref{full-cyclic},
let $ax_1=x_2b$ for some $a,b, x_1, x_2 \in X$. Then
\begin{enumerate}
\item  there exist
$c, x_3, x_0 \in X$ such that: (a) $ax_2=x_3b$, (b)  $ax_0=x_1b$,
and (c)
 $cx_1=x_2a$,
\item if
(i) $ax_1=x_2b$, (ii) $ax_2=x_3b \; $   and (iii) $cx_1=x_2a$,
then $cx_2=x_3a$.
\end{enumerate}
\end{lemma}

\noindent {\bf Proof.} 1(a) and (b)  follow immediately from
condition (C). Indeed, (C) applied to $ax_1=x_2b$ implies $a*
=x_1b$ and $ax_2 = *b$, with "$*$" meaning an element of $X$. In
general the letter $*$ is different in the first and the second
equality. For (c) consider $x_2b = ax_1$. Applying (a) we get
$x_2a = cx_1$ for some $c\in X$.

2) Assume (i), (ii) and (iii) hold. Applying 1(c) to $ax_2=x_3b$
yields $t\in X$ such that $t x_2= x_3a$. But then, applying 1(b),
we get $ts=x_2a$ for some $s \in X$. Since $S$ is of skew type,
comparing the latter with $cx_1=x_2a$, we obtain $t=c$. Hence
$cx_2=x_3a$. $\Box$

\vspace{10pt} Now the statement of the proposition can be derived
from the lemma as follows.

Let $x,y \in X$. From (C) it follows that the sequence for the
"internal cycle" $x=x_1,\ldots , x_k$ exists, so that
$yx_i=x_{i+1}z$ and $yx_{k}=x_{1}z$, for some $z \in X$ and all
$i=1,\ldots ,k-1$. By 1(c) there exists $y^{(1)}\in X$, such that
$$(iv) \hspace{40pt} y^{(1)}x_1=x_2y . \hspace{100pt}$$ Hence
$yx_{1}=x_{2}z, yx_{2}=x_{3}z$ and $y^{1}x_{1}=x_{2}y$ (if $k=1$
then we put $x_{2}=x_{3}=x_{1}$ and if $k=2$ then we put
$x_{3}=x_{1}$). So because of 2) we get $y^{(1)}x_{2}=x_{3}y$. It
follows by an induction procedure that $y^{(1)}$ is compatible
with the whole cycle $x_{1},\ldots ,x_{k}$, that is
 $$y^{(1)}x_1=x_2y,  y^{(1)}x_2=x_3y, \ldots,y^{(1)}x_k=x_1y.$$
Applying the same procedure to (iv) we obtain a $y^{(2)}\in X$
such that $$ y^{(2)}x_1=x_2y^{(1)},  y^{(2)}x_2=x_3y^{(1)}, \ldots
, y^{(2)}x_k=x_1y^{(1)}.$$ Condition (C) applied to $x_{2}y$
implies that, after finitely many such steps we shall close the
cycle for $y$'s, that is, we obtain a sequence of pairwise
distinct $y^{(1)}, \ldots , y^{(p-1)}$ such that
$x_2y^{(i)}=y^{(i+1)}x_1,$ for $i=1,\ldots ,p-2$ and
$x_{2}y^{(p-1)}=yx_{1}$. Since $yx_{1}=x_{2}z$, we get
$y^{(p-1)}=z$. Also $$ y^{(i+1)}x_1=x_2y^{(i)},
y^{(i+1)}x_2=x_3y^{(i)} , \ldots , y^{(i+1)}x_k=x_1y^{(i)} $$ for
$i=1,\ldots ,p-2$.  The assertion follows by reindexing the
elements $y^{(1)},\ldots $, $y^{(p-1)}$.
 $\Box$

\vspace{10pt} The following result allows to construct many
examples of noetherian PI algebras from semigroups of skew type.

\begin{proposition}          \label{cyclic}
Assume that $S=\langle x_{1},\ldots ,x_{n}\rangle$ is a semigroup
of skew type that satisfies the cyclic condition and $S=\{
x_{1}^{a_{1}}\cdots x_{n}^{a_{n}} | a_{i}\geq 0 \}$. Then $K[S]$
is a finite left and right module over a commutative subring of
the form $K[A]$, where $A=\langle x_{1}^{p},\ldots
,x_{n}^{p}\rangle$ for some $p\geq 1$. Namely $S=\bigcup _{c\in
C}cA$ with $C=\{ x_{1}^{i_{1}}\cdots x_{n}^{i_{n}}\, | \, i_{j}<p
\}$ and $cA=Ac$ for every $c\in C$. In particular $K[S]$ is a
right and left noetherian PI algebra.
\end{proposition}
\noindent {\bf Proof.} Let $x,y_{1}\in X=\{ x_{1},\ldots ,x_{n}
\}$. Then, for some $t,y_{1},\ldots ,y_{s}\in X$ we have
\begin{equation}
\begin{array}{c}
 xy_{1}=y_{2}t\\
 xy_{2}=y_{3}t \\
 \vdots \\
 xy_{s}=y_{1}t.
 \end{array}
 \end{equation} This easily implies that
$x^{s}y_{i}=y_{i}t^{s}$ for all $i=1,\ldots ,s$. Hence for every
$x,y\in X$ there exists $t\in S$ such that $x^{r}y=yt^{r}$, where
$r$ is the least common multiple of lengths of all cycles in $S$.
Suppose there exist $y',y''$ so that $xy=y't$ and $x'y=y''t$. By
the cyclic condition we get $xy'''=yt$ and $x'y''''=yt$ for some
$y''',y''''$. Since the relations are of skew type this yields
$x=x'$. So the generator $y$ acts as an injection, and thus as a
bijection on the set $X$, by mapping $x$ to $t$ if $xy=y't$ for
some $y'\in X$. It then also follows that there is a multiple $p$
of $r$ such that $x_{i}^{p}x_{j}^{p}=x_{j}^{p}x_{i}^{p}$ for all
$i,j$. Every $y\in X$ acts also as a bijection on the set
$\{x_{1}^{p}, \cdots ,x_{n}^{p}\}$. Since $S=\{
x_{1}^{a_{1}}\cdots x_{n}^{a_{n}} | a_{i}\geq 0 \}$, it now
follows that $S=CA$. Moreover $cA=Ac$ for $c\in C$ because $c$
acts as a bijection on the set of generators of $A$. $\Box$

\vspace{10pt} We note that the previous proof still works if $S$
is a semigroup of skew type that satisfies the cyclic condition
and $S$ is the union of sets of the form $\{ y_{1}^{a_{1}}\cdots
y_{k}^{a_{k}} | a_{i}\geq 0 \}$ where $y_{1},\ldots, y_{k}\in X$
and $k\leq n$. In Theorem~\ref{general} we will prove that the
latter is a consequence of the cyclic condition. Moreover (see
Theorem~\ref{noether}) $K[S]$ is still a noetherian PI algebra for
a class of semigroups of skew type essentially wider than those
satisfying the cyclic condition.

\section{Noetherian implies PI}

It is well known that the Gelfand-Kirillov dimension of a finitely
generated PI-algebra is finite (see \cite{krause}). One of our
aims is to show that the converse holds for every algebra $K[S]$
of a semigroup $S$ of skew type, provided that $K[S]$ is right
noetherian. Surprisingly, the following theorem shows that this
can be proved in the more general context of finitely generated
monoids  defined by homogeneous relations. Clearly, in such a
semigroup we have  a natural degree function given by $s\mapsto
|s|$, where $|s|$ is the length of $s\in S$ as a word in the
generators of $S$.

In the proof of the theorem we rely on the rich structure of
linear semigroups (\cite{book2}).

\begin{theorem}     \label{growth}
Let $S$ be a monoid such that the algebra $K[S]$ is right
noetherian and $\GK (K[S])< \infty $. Then $S$ is finitely
generated. If, moreover, $S$ has a monoid presentation of the form
$$S=\langle x_{1},\ldots ,x_{n} \, | \, R \rangle $$ with $R$ a
set of homogeneous relations, then $K[S]$ satisfies a polynomial
identity.
\end{theorem}
\noindent {\bf Proof.}  The first assertion follows from
Theorem~2.2 in \cite{noeth}. So assume $S$ has a monoid
presentation $S=\langle x_{1},\ldots ,x_{n} \, | \, R \rangle $.
Note that the unit group $U(S)$ is trivial. Let $T=S^{0}$, the
semigroup with zero $\theta$ adjoined. We define a congruence
$\rho$ on $T$ to be homogeneous if $s\rho t$ and $(s,\theta )\not
\in \rho $ imply that $|s|=|t|$.

The contracted semigroup algebra $K_{0}[T]$ may be identified with
$K[S]$.  Suppose that $K_{0}[T]$ is not a PI algebra. Then, by the
noetherian condition, there exists a maximal homogeneous
congruence $ \eta $ on $T$ such that $K_{0}[T/\! \eta ]$ is not
PI. So, replacing $T$ by $T/\! \eta $, we may assume that every
proper homogeneous homomorphic image of $T$ yields a PI algebra.

Since there are only finitely many minimal prime ideals of
$K_{0}[T]$ and the prime radical $B(K_{0}[T])$ is nilpotent, there
exists a minimal prime $P$ such that $K_{0}[T]/P$ is not a PI
algebra. As $K_{0}[T]$ can be considered in a natural way as a
${\bf Z}$-graded algebra (with respect to the length function on
$S$), it is well known \cite{graded}, that $P$ is a homogeneous
ideal of $K_{0}[T]$. Therefore the congruence $\rho _{P}$
determined by $P$ is homogeneous. (Recall that $s\rho _{P}t$ if
$s-t\in P$, for $s,t\in T$.) Since $K_{0}[T]/P$ is a homomorphic
image of $K_{0}[T/\rho _{P}]$, and because of the preceding
paragraph of the proof, we get that $T=T/\rho _{P}$. As $K_{0}[T]$
is right noetherian, we thus get $$T\subseteq K_{0}[T]/P \subseteq
M_{t}(D)$$ for some division algebra $D$, where
$M_{t}(D)=Q_{cl}(K_{0}[T]/P)$, the classical ring of quotients of
$K_{0}[T]/P$. Let $I$ be the set of all elements of $T$ (with
$\theta $) that are of minimal nonzero rank as matrices in
$M_{t}(D)$. Consider $K\{ I \}$, the subalgebra of $K_{0}[T]/P$
generated by $I$. Clearly $K\{ I \}$ is an ideal of $K_{0}[T]/P$.
Then $M_{t}(D)=Q_{cl}(K\{ I\})$. So $K\{ I\}$ is not a PI algebra,
as otherwise its ring of quotients would also satisfy a polynomial
identity.

Since not all elements of $I$ can be nilpotent, it follows from
the theory of linear semigroups that $I$ has
 a nonempty intersection $C$  with a maximal subgroup
$G$ of the multiplicative monoid $M_{t}(D)$. So $G$ is the group
of units of the monoid $eM_{t}(D)e$ for some $e=e^{2}$ in
$M_{t}(D)$.  Let $F\subseteq G$ be the group generated by $C$.
Define $$Z=\{ ex | x\in T, Cx\subseteq C \} .$$ If $g=ex\in G\cap
eT$ then $Cx=Cex\in G\cap T=C$. Hence $g\in Z$ and $G\cap
eT\subseteq Z$. It is easy to see that $Z\subseteq G$, so that
$Z=G\cap eT$. We claim that the monoid $Z$ satisfies the ascending
chain condition on right ideals. Fix some $c\in C$. Let $J$ be a
right ideal of $Z$. Notice that $cJT$ is a right ideal of $T$.
Then $cJT\cap Z=cJeT\cap Z=cJZ=cJ$ because $cJet\in Z$ implies
$et\in G\cap eT=Z$ for $t\in T$. As $T$ is a cancellative monoid
with the ascending chain condition on right ideals, the claim
follows.

One verifies that $F$ is a finitely generated group. This follows
from Proposition~3.16 in \cite{book2} (the result is proved for a
field $D$ only, but the proof works also for division rings $D$.)
Since $\GK (K_{0}[T])<\infty $ we also have $\GK (C)<\infty $. It
is then known that $F$ also has finite Gelfand-Kirillov dimension,
\cite{grigor}. Moreover, as $F$ is finitely generated, it follows
from \cite{gromov} that $F$ is nilpotent-by-finite.

Next we claim that the group of units $U(Z)$ of $Z$ is a periodic
group. For this, suppose $g, g^{-1}\in Z$. Then $Cg\subseteq C$
and $Cg^{-1}\subseteq C$. So $Cg=C$. Write $g=ab^{-1}$ with
$a,b\in C$. Then $Ca=Cb$ and so $Ma=Mb$, where $M$ is the subset
consisting of the elements of minimal length in $C$. Clearly
$Mg=M$. As $M$ is finite, we get $g^{k}=e$ for some $k\geq 1$,
which proves the claim.

So $U(Z)$ is a periodic subgroup of the finitely generated
nilpotent-by-finite group $F$. Hence $U(Z)$ is finite. Since also
$Z$ satisfies the ascending chain condition on right ideals, it
follows from the remark on page 550 in \cite{noeth} that $F$ is
finite-by-abelian-by-finite. Hence $F$ is abelian-by-finite and
thus $K[F]$ is a PI algebra.

Finally, as $T$ satisfies the ascending chain condition on right
ideals, $I$ intersects finitely many ${\cal R}$-classes of the
monoid $M_{t}(D)$. It is then known that $I$ embeds into a
completely $0$-simple semigroup with finitely many ${\cal
R}$-classes and with a maximal subgroup $F$. It  follows that $K\{
I \}$ is a PI algebra, see \cite{book}, Proposition~20.6, a
contradiction. This completes the proof of the theorem. $\Box $

\begin{corollary}         \label{pi}
Let $S$ be a semigroup of skew type such that $\GK (K[S])<\infty
$.
 If $K[S]$ is right noetherian, then it
satisfies a polynomial identity. In particular, $K[S]$ embeds into
a matrix ring over a field and $\GK (K[S])=\GK (K[S]/B(K[S]))$ is
an integer, where $B(K[S])$ is the prime radical of $K[S]$.
Moreover, $S$ satisfies a semigroup identity.
\end{corollary} \noindent {\bf Proof.} $K[S]$ is a PI algebra by
Theorem~\ref{growth}. Hence \cite{ananin} implies that $K[S]$ is a
subalgebra of $M_{t}(L)$ for a field $L$ and $t\geq 1$. Then, by a
result of Markov $\GK (K[S])=\GK (K[S]/B(K[S]))$ is an integer,
see \cite{krause}, Section~12.10. The last assertion now follows
from \cite{book2}, Proposition~7.10. $\Box$

\vspace{10pt} In Sections~\ref{sec4} and \ref{sec5} we will show
that for a wide class semigroups $S$ of skew type $\GK
(K[S])<\infty$ and $K[S]$ is right and left noetherian. So the
corollary is applicable in this situation.

\section{Non-degenerate and the ascending chain condition} \label{sec4}

Assume $S=\langle x_{1},\ldots ,x_{n}\rangle $ is a semigroup of
skew type that satisfies the cyclic condition. If $x\in \{
x_{1},\ldots ,x_{n} \}$, then for every $y_{1}\in \{ x_{1},\ldots
,x_{n} \}$ we get a cycle
 \begin{equation}
\begin{array}{c}
 xy_{1}=y_{2}t\\
 xy_{2}=y_{3}t \\
 \vdots \\
 xy_{s}=y_{1}t
 \end{array}
 \end{equation}
with $t,y_{i}\in \{ x_{1},\ldots ,x_{n} \}$. Since every $xx_{k}$,
with $x\neq x_{k}$, appears in one of the relations defining $S$,
it is clear that for every $x_{k}$ there exists a relation of the
form $xx_{i}=x_{k}x_{l}$ for some $i,l$.

A semigroup of skew type satisfying the latter condition will be
said right non-degenerate. Left non-degenerate semigroups are
defined dually. A symmetric argument shows that the cyclic
condition implies that $S$ is left non-degenerate as well. Notice
that if $S$ is right non-degenerate then every $x\in \{
x_{1},\ldots ,x_{n}\}$ defines a bijection $f_{x}$ of $\{
x_{1},\ldots ,x_{n} \}$ as follows: if $xx_{i}=x_{k}x_{l}$ then
$f_{x}(x_{i})=x_{k}$.

There are many examples of right and left non-degenerate $S$ which
do not satisfy the cyclic condition. For example, $S=\langle
x_{1},x_{2},x_{3},x_{4}\rangle $ defined by the relations:
$x_{2}x_{1}=x_{1}x_{3},x_{3}x_{1}=x_{2}x_{4},x_{4}x_{1}=x_{1}x_{2},
x_{3}x_{2}=x_{1}x_{4},x_{4}x_{2}=x_{2}x_{3},x_{4}x_{3}=x_{3}x_{4}$.

First we prove some technical and combinatorial properties of
non-de\-gen\-erate semigroups.

\vspace{10pt} Let $S$ be a semigroup of skew type. Let $Y=\langle
X \rangle, X= \{ x_{1},\ldots ,x_{n}\} $, be a free monoid of rank
$n$. (So, we use the same notation for the generators of $Y$ and
of $S$, if unambiguous.) For any $m\geq 2$ and any $y_{1},\ldots ,
y_{m}\in X$ define $$g_{i}(y_{1}\cdots y_{m})=y_{1}\cdots
y_{i-1}\overline{y}_{i}\overline{y}_{i+1}y_{i+2}\cdots y_{m}$$ for
$i=1,\ldots , m-1$, where
$$y_{i}y_{i+1}=\overline{y}_{i}\overline{y}_{i+1}$$ is one of the
defining relations of $S$ (if $y_{i}\neq y_{i+1}$) or
$y_{i}=y_{i+1}=\overline{y}_{i}=\overline{y}_{i+1}$. Let
$$g(y_{1}\cdots y_{m})=g_{m-1}\cdots g_{2}g_{1}(y_{1}\cdots y_{m})
$$ (Notice that $g$ is used for all $m=2,3,\ldots $.) If
$g(y_{1}\cdots y_{m})=s_{1}\cdots s_{m}, \, s_{i}\in X$, then we
set $$f_{y_{1}}(y_{2}\cdots y_{m})=s_{1}\cdots s_{m-1} .$$ So
$f_{y_{1}}:X^{m-1}\longrightarrow X^{m-1}$ can be considered as a
function on the subset $X^{m-1}$ of $Y$ consisting of all words of
length $m-1$.

\begin{lemma}      \label{one-to-one}
Assume that $S$ is a right non-degenerate semigroup of skew type.
If $y_{1}\in X$, then $f_{y_{1}}:X^{m-1}\longrightarrow X^{m-1}$
is a one-to-one mapping, for any $m\geq 2$.
\end{lemma}
\noindent  {\bf Proof.} We proceed by induction on $m$. The case
$m=2$ is clear because $S$ is right non-degenerate.

Assume now that $m>2$. Let $$g(y_{1}\cdots y_{m})=s_{1}\cdots
s_{m}. $$ We will show that $s_{1}\cdots s_{m-1}$ and $y_{1}$
determine $y_{2}\cdots y_{m}$. Notice that
\newline $g_{1}(y_{1}\cdots y_{m})=s_{1}h(y_{1}y_{2})y_{3}\cdots
y_{m}$ where $$y_{1}y_{2}=s_{1}h(y_{1}y_{2})$$ is a relation in
$S$ or $y_{1}=y_{2}=s_{1}=h(y_{1}y_{2})$. Moreover
$$s_{1}g(h(y_{1}y_{2})y_{3}\cdots y_{m})=g(y_{1}\cdots
y_{m})=s_{1}\cdots s_{m} .$$ Then $g(h(y_{1}y_{2})y_{3}\cdots
y_{m})=s_{2}\cdots s_{m}$ and hence by the induction hypothesis it
follows that $s_{2}\cdots s_{m-1}$ and $h(y_{1}y_{2})$ determine
$y_{3}\cdots y_{m}$. Since $S$ is right non-degenerate, $y_{1}$
and $s_{1}$ determine $h(y_{1}y_{2})$ and $y_{2}$. Hence $y_{1}$
and $f_{y_{1}}(y_{2}\cdots y_{m})=s_{1}\cdots s_{m-1}$ determine
$y_{2}\cdots y_{m}$, as desired. $\Box$

\vspace{10pt} Our aim is to investigate when $K[S]$ is noetherian.
Hence we first study the weaker condition that $S$ satisfies the
ascending chain condition on right ideals.

Let $S=\langle x_{1},\ldots ,x_{n}\rangle$. We shall consider the
following over-jumping property
\begin{quote}
for every $a\in S$ and every $i$ there exist $k\geq 1$ and $w\in
S$ such that $aw=x_{i}^{k}a$.
\end{quote}
This property is formally stronger than the following immediate
consequence of the ascending chain condition on right ideals in
$S$
\begin{quote}
for every $a\in S$ and every $i$ there exist positive integers
$q,p$ and $w\in S$ such that $x_{i}^{p}aw=x_{i}^{p+q}a$.
\end{quote}
(Indeed, this condition immediately follows from the ascending
chain condition applied to $I_{j}=\bigcup _{k=1}^{j}x_{i}^{k}aS$.)
We show that the over-jumping property holds for the class of
right non-degenerate semigroups of skew type.

\begin{proposition}        \label{over}
Assume that $S$ is a right non-degenerate semigroup of skew type.
Then $S$ has the over-jumping property.
\end{proposition}
{\bf Proof.} Fix some $y_{1}\in X=\{x_{1},\ldots ,x_{n} \}$. We
have shown that, if $m\geq 2$, then
$f=f_{y_{1}}:X^{m-1}\longrightarrow X^{m-1}$ is a permutation.
Therefore $f^{r}$ is the identity map for some $r\leq
(|X|^{m-1})!=(n^{m-1})!$. So, for any $y_{2},\ldots ,y_{m}\in X$
we have $$f^{r}(y_{2}\cdots y_{m})=y_{2}\cdots y_{m}.$$ Now,
interpreting $X$ as the generating set of $S$, we get the
following equality in $S$ $$y_{1}y_{2}\cdots
y_{m}=f_{y_{1}}(y_{2}\cdots y_{m})s_{m} .$$ Next
$$y_{1}^{2}y_{2}\cdots y_{m}=y_{1}f_{y_{1}}(y_{2}\cdots
y_{m})s_{m}=f_{y_{1}}(f_{y_{1}}(y_{2}\cdots y_{m}))s_{m+1}s_{m}$$
for some $s_{m+1}\in X$. Proceeding this way, we come to
$$y_{1}^{r}y_{2}\cdots y_{m}=f_{y_{1}}^{r}(y_{2}\cdots
y_{m})s_{m+r-1}\cdots s_{m+1}s_{m}=y_{2}\cdots
y_{m}s_{m+r-1}\cdots s_{m+1}s_{m}$$ for some $s_{i}\in X,
i=m,\ldots ,m+r-1$. This means that in $S$ we have
$$y_{1}^{r}y_{2}\cdots y_{m}=y_{2}\cdots y_{m}w$$ for some $w\in
S$.

This can be also repeated for $f$ considered as a map $f_{y_{1}}:
X\cup \cdots \cup X^{m-1}\longrightarrow  X\cup \cdots \cup
X^{m-1}$. We have thus shown that $S$ has the following property:
\begin{quote}
 for every $m\geq 1$ there exists $r\geq 1$ ($r\leq (n^{m-1})!$)
 such that if $a\in S$
has length less than $m$ in the generators $x_{1}, \ldots ,x_{n}$
and $i\in \{ 1,\ldots ,n \}$ then we have $aw=x_{i}^{r}a$ for some
$w\in S$.
\end{quote}
The result follows. $\Box $

\begin{lemma}          \label{nondegener}
Assume that $S$ is a right non-degenerate semigroup of skew type.
Then for every $x,y\in S$ there exist $t,w\in S$ such that
$|w|=|y|$ and $xw=yt$.
\end{lemma}
\noindent {\bf Proof.}  Suppose first $|x|=1$, so that $x=x_{j}$
for some $j$. Then the assertion follows from
Lemma~\ref{one-to-one}. So, suppose $|x|>1$. We now proceed by
induction on the length of $x$ as a word in $x_{1},\ldots ,x_{n}$.
So suppose that the assertion holds for all $x\in S$ of length
$<m$. Let $x\in S$ be such that $|x|=m$, say $x=z_{1}\cdots z_{m}$
for some $z_{i}\in \{ x_{1},\ldots ,x_{n}\}$. By the induction
hypothesis $z_{1}\cdots z_{m-1}u=yw$ for some $u,w\in S$ with
$|u|=|y|$. We know also that $z_{m}v=us$ for some $v,s\in S$ such
that $|v|=|u|$. Then $$xv=z_{1}\cdots z_{m-1}z_{m}v=z_{1}\cdots
z_{m-1}us=yws.$$ Since $|v|=|y|$, this proves the assertion.
 $\Box $

\vspace{10pt} The following result, together with its proof,
provide the first insight into the structure of non-degenerate
semigroups and their algebras. This will be heavily exploited and
strengthened in Section~5.

In the proof the following sets will play a crucial role.

\begin{definition} \label{spset}
Let $S=\langle x_{1},x_{2}, \ldots , x_{n} \rangle$ be a semigroup
of skew type. For a subset $Y$ of $X=\{ x_{1}, \ldots , x_{n}\}$
define
 $$S_{Y}=\bigcap _{y\in Y} yS$$
and $$D_{Y}=\{ s\in S_{Y} \mid \mbox{ if } s=xt \mbox { for some }
  x\in X \mbox{ and } t\in S \mbox{ then } x\in Y\} .$$

The left-right symmetric duals of these sets will be denoted by
$S_{Y}'$ and $D_{Y}'$ respectively.
\end{definition}

Notice that because of Lemma~\ref{nondegener} each such set
$S_{Y}$ is non-empty. However, it may happen that  $S_{Y}=S_{Z}$
for different subsets $Y$ and $Z$ of $X$; possibly it can occur
that $D_{Y} = \emptyset$.

\begin{theorem}        \label{general}
Let $S=\langle x_{1},\ldots ,x_{n}\rangle $ be a semigroup of skew
type. If $S$ is right non-degenerate then
\begin{enumerate}
\item for each integer $i$, with $1\leq i \leq n$,
  $S_{i}=\bigcup_{Y : |Y|=i}S_{Y}$
is an ideal of $S$, and
 $$S_{X}=S_{n} \subseteq S_{n-1} \subseteq \cdots \subseteq S_{1}
 \subseteq S,$$
\item $S$ is the union of sets
of the form $\{y_{1}^{a_{1}}\cdots y_{k}^{a_{k}}:a_{i}\geq 0\}$,
where $y_{1},\ldots ,y_{k}$ $\in X$ and $k\leq n$.
\end{enumerate}
 In particular, $\GK (K[S])\leq n $.
\end{theorem}
{\bf Proof.} Let $Y$ be a subset of $X=\{ x_{1}, \ldots ,
x_{n}\}$.
 If $x\in X$ then let
$Z\subseteq X$ be the largest subset such that $xS_{Y}\subseteq
S_{Z}$. Since $X$ is right non-degenerate, it follows that
$|Z|\geq |Y|$. Moreover, if $x\not \in Y$, then $|Z|>|Y|$.
Consequently, $S_{j}=\bigcup_{Y : |Y|=j} S_{Y}$ are ideals of $S$
such that $$S_{X}=S_{n}\subseteq S_{n-1}\subseteq \cdots \subseteq
S_{1} \subseteq S.$$ Note that  if $j=|Y|$ then
$D_{Y}=S_{Y}\setminus S_{j+1}$ (we let $S_{n+1}=\emptyset $). So
$$S_{j}\setminus S_{j+1} = \bigcup _{Z : |Z|=j} D_{Z}$$ is a
disjoint union.

Suppose first that $|Y|=1$. Let $w\in D_{Y}\setminus \langle
y\rangle$, where $Y=\{ y\}$. Then $w=y^{k}xt$ for some $x\in X,
t\in S$, and $y^{k}x\in D_{Y}\setminus \langle y\rangle $. Since
$S$ is right non-degenerate, there exist $r\geq 1$ and distinct
elements $u_{1},\cdots ,u_{r}=x\in X$ such that $yx=u_{1}w_{1},
yu_{1}=u_{2}w_{2}, \ldots, yu_{r-1}=u_{r}w_{r}$. Therefore
$y^{q}x\in u_{1}S\cup \cdots \cup u_{r}S$ for every $q\geq 1$. But
$u_{i}\neq y$ for all $i\geq 1$, so $y^{k}x\not \in D_{Y}$, a
contradiction. It follows that $D_{Y}=\langle y \rangle \setminus
\{ 1\} $.

Fix some $y\in Y$. Suppose $s\in D_{Y}$ and $j=|Y|$. Let $r\geq 1$
be the maximal integer such that $s=y^{r}t$ for some $t\in S$.
Suppose $t\in S_{Z}$ for some $Z\subseteq X$ with $|Z|=|Y|$. If
$y\not \in Z$ then $yt\in x_{i}S$ for at least $|Y|+1$ different
indices $i$. So $yt\in S_{j+1}$ and therefore $s\in S_{j+1}$, a
contradiction. So, we have $y\in Z$. Then $t\in yS$, which
contradicts the maximality of $r$. Hence, we have shown that
$t\not \in S_{j}$. It follows that $$D_{Y}\subseteq \langle y
\rangle (S\setminus S_{j}).$$ By induction on $|Y|$ this easily
implies $D_{Y}$ is contained in a union of sets of the form $\{
y_{1}^{a_{1}}\cdots y_{j}^{a_{j}}: a_{j}\geq 0\}$, where $|Y|=j$
and $y_{i}\in X$. So $S$ is the (finite) union of sets of the form
$\{y_{1}^{a_{1}}\cdots y_{k}^{a_{k}}:a_{i}\geq 0\}$, where
$y_{1},\ldots ,y_{k}\in X$ and $k\leq n$.

The assertion on the Gelfand-Kirillov dimension of $K[S]$ is now
an easy consequence. It is clear that $S\setminus S_{2}=\bigcup
_{i=1}^{n}\langle x_{i}\rangle $. Hence there are $nm+1$ elements
of $S$ that are words of length at most $m$ in the generators
$x_{1},\ldots ,x_{n}$ and that lie in $S\setminus S_{2}$.
Proceeding by induction on $j$, assume that the number of elements
of $S\setminus S_{j}$ that are words of length at most $m$ is
bounded by a polynomial of degree $j-1$ in $m$. Let $|Y|=j,
Y\subseteq X$. Since $D_{Y}\subseteq \langle y \rangle (S\setminus
S_{j})$ for $y\in Y$, it is easy to see that the number of
elements of $D_{Y}$ that are words of length at most $m$ is
bounded by a polynomial of degree $j$. As $S_{j}\setminus S_{j+1}$
is a finite union of such $D_{Y}$, the same is true of the
elements of the set $S_{j}\setminus S_{j+1}$. This proves the
inductive claim. It follows that the growth of $S$ is polynomial
of degree not exceeding $n$, so that $\GK (K[S])\leq n$. $\Box$

\vspace{10pt} The left-right symmetric dual of $S_{i}$ will be
denoted by $S_{i}'$. Of course, if $S$ is a semigroup of skew type
which is left non-degenerate then we obtain that each $S_{i}'$
also is an ideal of $S$.

\vspace{10pt}
 The following technical
result turns out to be very useful.

\begin{lemma}           \label{claim}
Let $S$ be a right non-degenerate semigroup of skew type. Let $Y$
be a subset of $X$ and assume $|Y|=i-1$. Let $b\in D_{Z}$, for
some subset $Z$ of $Y$. Assume that $k$ is the length of $b$ in
the generators of $X$. Then $(S_{i-1})^{k}\cap D_{Y} \subseteq
bS$. Furthermore, $(S_{i-1})^{k+1}\cap D_{Y} \subseteq bS_{i-1}$
and $(S_{i-1}\cap S_{i-1}')^{k+1} \cap D_{Y} \subseteq b (S_{i-1}
\cap S_{i-1}')$.
\end{lemma}
{\bf Proof.} If $k=1$ the assertion is clear. So assume $k\geq 2$.
Write $b=y_{k}\cdots y_{1}$ with each $y_{j}\in X$. Let $q\geq k$
and  $a=a_{q}\cdots a_{1}\in D_{Y}$ with each $a_{j}\in
S_{i-1}\setminus S_{i}$. Since $b\in D_{Z}$ with $Z\subseteq Y$ we
get $y_{k}\in Y$, and therefore $a_{q}\in y_{k}S$. So
$a_{q}=y_{k}b_{k}$ for some $b_{k}\in S$. Then
$a_{q}a_{q-1}=y_{k}c_{k}$ where $c_{k}=b_{k}a_{q-1}$. Clearly
$c_{k}\in S_{i-1}\setminus S_{i}$.

Suppose we have already shown that
   \begin{eqnarray}
     a_{q}\cdots a_{q-r} &=&y_{k}\cdots y_{k-r+1}c_{k-r+1}
     \label{moving}
   \end{eqnarray}
for some $r\geq 1$ and $c_{k-r+1}\in S_{i-1}\setminus S_{i}$. We
claim that $c_{k-r+1}\in y_{k-r}S$. Let $W\subseteq X$ be so that
$|W|=i-1$ and $c_{k-r+1}\in D_{W}$. Consider the set
 $$U=\{x \in X \mid y_{k}\cdots y_{k-r+1}x\in D_{V} \mbox{ for some }
  V\subseteq Y\} .$$
Because of the right non-degeneracy, an induction argument on $r$
yields that $|U|\leq |Y|$. Since the left hand side of equation
(\ref{moving}) is an initial segment of $a$ and $a\in D_{Y}$ it
follows that $W\subseteq U$. So $W=U$. Since $y_{k} \cdots
y_{k-r+1}y_{k-r}$ is an initial segment of $b\in D_{Z}$ and
$Z\subseteq Y$ we also get that $y_{k-r}\in U=W$. Hence
$c_{k-r+1}\in y_{k-r}S$. This proves the claim.

Now write $c_{k-r+1}=y_{k-r}b_{k-r}$ for some $b_{k-r}\in S$. So
\begin{eqnarray*}
     a_{q}\cdots a_{q-r}a_{q-r-1} &=&y_{k}\cdots y_{k-r+1}y_{k-r}b_{k-r}a_{q-r-1}.
   \end{eqnarray*}
Define $c_{k-r}=b_{k-r}a_{q-r-1}$. Then $c_{k-r}\in
S_{i-1}\setminus S_{i}$.

So we have shown that for any $q\geq k$, $a_{q}\cdots a_{q-k+1}
\in y_{k}\cdots y_{1}S=bS$. If $q=k$ then the first assertion of
the lemma follows. On the other hand, if $q=k+1$ then we obtain
$a=a_{k+1}\cdots a_{1} \in bSa_{1} \subseteq bS_{i-1}$. The second
and third assertion of the lemma now easily follow. $\Box$

\begin{proposition}         \label{acc}
Let $S$ be a right non-degenerate semigroup of skew type. Then $S$
has the ascending chain condition on right ideals.
\end{proposition}
{\bf Proof.} Suppose we know already that $S/S_{i}$ has the
ascending chain condition on right ideals for some $i$. We will
show that $S/S_{i+1}$ also has this property. Recall that by
definition $S_{n+1}=\emptyset $ and $S/S_{n+1}=S$. Then with
$i=n+1$ the assertion follows.

  From Theorem~\ref{general} we know that $S=\{
z_{1}^{a_{1}}z_{2}^{a_{2}}\cdots z_{m}^{a_{m}} : a_{j}\geq 0
 \}$ for some $m\geq 1$ and $z_{1},\ldots ,z_{m}\in X$
(not all $z_{j}$ are necessarily different). We claim that
$S_{i}/S_{i+1}$ is a finitely generated as a right ideal of
$S/S_{i+1}$. To prove this, it is sufficient to show by induction
on $m-k$ that right ideal of $S/S_{i+1}$ generated by $C_{k}\cap
(S_{i}\setminus S_{i+1})$ is finitely generated; where $C_{k}=\{
z_{k}^{a_{k}}\cdots z_{m}^{a_{m}} | a_{j}\geq 0 \}$. The case
$m-k=0$ is clear. The case $m-k=m-1$ gives the assertion.

So assume $1\leq k< m$. Let $B=\{ b\in C_{k+1} | z_{k}^{a}b \in
S_{i} \mbox{ for some } a \}$. If $y\in C_{k}\cap (S_{i}\setminus
S_{i+1})$ then $y\in z_{k}^{a}(B\cap (S\setminus S_{i}))$ for some
$a$ (see the proof of Theorem 4.1) and $(B\cap (S\setminus
S_{i}))S\subseteq b_{1}S\cup \cdots \cup b_{r}S$ for some
$b_{j}\in B\cap (S\setminus S_{i})$ because $S/S_{i}$ has the
ascending chain condition on right ideals. Since $S_{i}$ is an
ideal of $S$, it follows that
\begin{equation}     \label{union}
C_{k}\cap (S_{i}\setminus S_{i+1}) \subseteq \bigcup _{t\geq N}
\bigcup _{j=1}^{r} z_{k}^{t}b_{j}S \cup \bigcup _{j=0}^{N-1}
z_{k}^{j}B_{j}
\end{equation}
where $N$ is chosen so that $z_{k}^{N}b_{j}\in S_{i}$ for
$j=1,\ldots ,r$ and $B_{j}=\{ y\in C_{k+1} | z_{k}^{j}y\in S_{i}
\}$. By the inductive hypothesis every $B_{j}\cap S_{i}$ generates
a finitely generated right ideal modulo $S_{i+1}$. On the other
hand, $(B_{j}\setminus S_{i})S$ is a finitely generated right
ideal because $S/S_{i}$ has the ascending chain condition on right
ideals. Hence $B_{j}$ and thus also $z_{k}^{j}B_{j}$ generates a
finitely generated right ideal modulo $S_{i+1}$.

Next we show that the double union above is a finitely generated
right ideal of $S$. Because of Proposition~\ref{over} we know that
$S$ has the over-jumping property. Consequently, for every $j$
there exist $w_{j}\in S$ and a positive integer $q_{j}$ such that
$$b_{j}w_{j}=z_{k}^{q_{j}}b_{j}.$$ Hence $$z_{k}^{q_{j}{\bf
N}}b_{j}\subseteq b_{j}S \, \mbox{ and so } \, z_{k}^{p+q_{j}{\bf
N}}b_{j}\subseteq z_{k}^{p}b_{j}S $$ for every $p\geq 0$. It
follows that the right ideal $\bigcup _{t\geq N}z_{k}^{t}b_{j}S=
\bigcup _{t=N}^{N+q_{j}} z_{k}^{t}b_{j}S$ is finitely generated,
as claimed.

As the left and the right side in (\ref{union}) generate modulo
$S_{i+1}$ the same right ideal, it follows that $C_{k}\cap
(S_{i}\setminus S_{i+1})$ generates a finitely generated right
ideal modulo $S_{i+1}$. So we proved our claim that
$S_{i}/S_{i+1}$ is a finitely generated right ideal of
$S/S_{i+1}$.

Suppose there is an infinite sequence $a_{1},a_{2},\ldots \in
S\setminus S_{i+1}$ such that we have proper inclusions
$$a_{1}S\subset a_{1}S\cup a_{2}S\subset \cdots \subset a_{1}S\cup
\cdots \cup a_{k}S\subset \cdots .$$ Since $S$ is the union of
finitely many sets $D_{Y}, Y\subseteq X$, we may assume that all
$a_{j}\in D_{Y}$ for some $Y$. As $S/S_{i}$ has the ascending
chain condition on right ideals, it follows that $D_{Y}\subseteq
S_{i}\setminus S_{i+1}$. Lemma~\ref{claim} implies that $a_{j}\not
\in S_{i}^{t}$ where $t$ denotes the length of $a_{1}$. This leads
to a contradiction with the fact that $S_{i}/S_{i+1}$ is a
finitely generated right ideal of $S/S_{i+1}$ and $S/S_{i}$ has
the ascending chain condition on right ideals. (Namely, if $S_{i}=
s_{1}S\cup \cdots \cup s_{q}S \cup S_{i+1}$ then $\{a_{j} \}$ has
a subsequence contained in $s_{k_{1}}\cdots s_{k_{p}}(S\setminus
S_{i})$ for some $p < t$ and some $k_{j}$, leading to a
contradiction.)

This proves that $S/S_{i+1}$ has the ascending chain condition on
right ideals, completing the inductive argument, and proving the
result. $\Box $

\section{Non-degenerate implies noetherian} \label{sec5}

Our main aim in this section is to show that left and right
non-degenerate semigroups $S$ yield left and right noetherian
algebras $K[S]$. To prove this we will rely on a general result
\cite{search} that makes use of ideal chains in $S$ of a special
type. Before stating the latter we recall some terminology. Let
$E={\cal M}(G,t,t;Id)$ be an inverse semigroup over a group $G$
with $t\geq 1$ (see \cite{howie}). In other words, $E=\{ (g)_{ij}
\mid g\in G , \; i,j=1, \ldots ,t \} \cup \{ 0 \}$, where
$(g)_{ij}$ denotes the $t\times t$-matrix with $g$ in the
$(i,j)$-component and zeros elsewhere. The multiplication on $E$
is the ordinary matrix multiplication. A semigroup $S$ is said to
be a generalised matrix semigroup if it is a subsemigroup of a
semigroup  $E$ of the above type and for every $i,j$ there exists
$g\in G$ so that $(g)_{ij}\in S$. So, in the terminology of
\cite{book}, $S$ is a uniform subsemigroup of $E$.

\begin{theorem}[Theorem 3.3 \cite{search}]     \label{sufficient}
Assume that $M$ is a finitely generated \newline monoid with an
ideal chain $M_{1}\subseteq M_{2}\subseteq \cdots \subseteq
M_{r}=M$ such that $M_{1}$ and every factor $M_{j}/M_{j-1}$ is
either nilpotent or a generalised matrix semigroup. If $M$ has the
ascending chain condition on right ideals, and $\GK (K[M])$ is
finite, then $K[M]$ is right noetherian.
\end{theorem}

 Now we assume that $S$ is a semigroup of skew type that
is right and left non-degenerate. Recall that for a
 $Y\subseteq X$ we
write $S_{Y}'= \bigcap _{y\in Y}Sy$. Further, $S_{k}'= \bigcup
_{Y\subseteq X, |Y|=k} S_{Y}'$ is an ideal of $S$. We claim that
$S_{k}^{k}\subseteq S_{k}'$ for all $k$ and symmetrically
$(S_{k}')^{k}\subseteq S_{k}$. Since $S_{1}=S_{1}'$ (see the
description of $S_{1}$ in the proof of Theorem~\ref{general}), we
may assume that $k>1$. Let $1\neq a\in S$ and $b\in S_{k}$. Let
$Y\subseteq X$ be maximal such that $a\in S_{Y}'$. If $|Y|<k$ then
there exists $x\in X$ such that $b\in xS$ and $x\not \in Y$. So we
may write $b=xc, c\in S$. Now $ax\subseteq Sx$ but, as $x\not \in
Y$ and $a\in S_{Y}'$, the element $ax$ is also contained in $|Y|$
different left ideals of the form $Sw$, $x\neq w\in X$ (use the
left non-degeneracy of $S$). Therefore $ax\in S_{Z}'$ for some
$Z\subseteq X$ with $|Z|>|Y|$. Hence $ab=axc\in S_{|Z|}'$. By
induction it follows easily that $S_{k}^{k}\subseteq S_{k}'$, as
desired.

We have shown that $I_{k}=S_{k}\cap S_{k}'$ is an ideal of $S$
such that $S_{k}/I_{k}$ is nilpotent.

\begin{theorem}        \label{noether}
Let $S=\langle x_{1},\ldots ,x_{n}\rangle $ be a semigroup of skew
type. If $S$ is right and left non-degenerate, then $K[S]$ is a
right and left noetherian PI-algebra.
\end{theorem}
{\bf Proof.}
 From Theorem~\ref{general} we know that $\GK (K[S])$
is finite. Because of Proposition~\ref{acc} we also know that $S$
satisfies the ascending chain condition on one sided ideals. In
view of Theorem~\ref{sufficient} and its dual, to prove that
$K[S]$ is right and left noetherian it is sufficient to show that
$S$ has an ideal chain with each factor either nilpotent or a
generalised matrix semigroup.

Write $S_{n+1}=S_{n+1}'=\emptyset$ and adopt the convention
$S/\emptyset=S$. By induction on $i$ we will prove that $S/S_{i}$
has an ideal chain of the desired type. The case $i=n+1$ then
yields the result. As noticed in the proof of
Theorem~\ref{general}, $S\setminus S_{2}$ is the disjoint union of
all $D_{\{ x_{i}\} } =\langle x_{i}\rangle \setminus \{ 1 \}$ and
$\{ 1 \}$. So $S/S_{2}$ has an ideal chain with commutative
$0$-cancellative factors, hence it has a chain of the type
described in Theorem~\ref{sufficient}. So now assume that we have
shown this for the semigroup $S/S_{i-1}$ for some $i\geq 3$.

Let $J'$ be the ideal of $S$ such that $S_{i}\subseteq J'$ and
$J'/S_{i}$ is the maximal nil ideal of $S/S_{i}$. We consider the
following ideals of $S$
 $$
  S_{i}\cap I \subseteq (S_{i}\cup S_{i}')\cap I \subseteq J
  \subseteq I=S_{i-1}\cap S_{i-1}'\subseteq S_{i-1}
 $$
where $J=J'\cap I$. (Notice that the first and the last Rees
factor are nilpotent by the comment after
Theorem~\ref{sufficient}.) Then $J/S_{i}$ is nilpotent because of
the ascending chain condition on one-sided ideals in $S$, see
Theorem~17.22 in \cite{faith}.

For $Y\subseteq X$, let  $D_{Y}$ and $D_{Y}'$
 be the subsets of $S$ introduced in Definition~\ref{spset}.
Let $I_{YW}=D_{Y}\cap D_{W}'$ for $Y,W\subseteq X$ and
$I_{Y}=I_{YY}$. If $|Y|=i-1$, $D_{Y}' \neq \emptyset$ and $x\in
X\setminus Y$ then $D_{Y}'x\subseteq S_{i}'$. (Use the left
non-degeneracy of $S$; $D_{Y}'x\subseteq Sx$ but also it is
contained in $|Y|$ different left ideals of the form $Syx$, $y\in
Y$, and thus $yx=x'y'$ with $y'\neq x$.) So $D_{Y}'D_{Z}\subseteq
S_{i}'$ for every $Y,Z\subseteq X$ of cardinality $i-1$ with
$Z\neq Y$, provided that $D_{Y}\neq \emptyset$ and $D_{Z}'\neq
\emptyset$. Hence we get a generalised matrix structure
$(S_{i-1}\cap S_{i-1}')\setminus (S_{i}\cup S_{i}') =\bigcup
_{Y,Z}I_{YZ}$.

Now, for $Y\subseteq X$ with $|Y|=i-1$ there are two mutually
exclusive cases:

\begin{enumerate}
 \item[Case]1: either $I_{Y}=\emptyset$ or there exist $b\in D_{Y}'$
  and $x\in D_{Y}$ such that $bx\in J$.
 \item[Case]2: $I_{Y} \neq \emptyset$ and
      $D_{Y}'D_{Y}\subseteq (S_{i-1}\cap S_{i-1}')\setminus J$,
      so in particular $I_{Y}$ is a subsemigroup
      of $(S_{i-1}\cap S_{i-1}')\setminus J$.
\end{enumerate}
In Case 1 we claim that $D_{Y}\cap I$ and $D_{Y}'\cap I$ are
contained in $J$. If $I_{Y}=\emptyset$, then the generalised
matrix structure easily yields that $(D_{Y}\cap I)^{2}$ and
$(D_{Y}'\cap I)^{2}$ are contained in $J$. As both $D_{Y}\cap I$
and $D_{Y}'\cap I$ are one sided ideals modulo $J$, we get
$D_{Y}\cap I ,\; D_{Y}'\cap I \subseteq J$. So assume $I_{Y}\neq
\emptyset$ and that there exist $b\in D_{Y}'$ and $x\in D_{Y}$ so
that $bx\in J$. Let $q$ be the maximum of the lengths of $b$ and
$x$. Then by Lemma~\ref{claim} and its right-left dual we get
$$I_{Y}^{2q}\subseteq (D_{Y}')^{q}D_{Y}^{q}\subseteq (Sb)(xS)\cup
S_{i} \cup S_{i}'\subseteq J'.$$ So $I_{Y}$ is nilpotent modulo
$J$. Then again $D_{Y}'\cap I$ (with zero) is a left ideal of
$S/S_{i}'$ and it is nil modulo $J$ (use the generalised matrix
pattern), so we must have $D_{Y}'\cap I \subseteq J$. Similarly
$D_{Y}\cap I \subseteq J$.

In Case~2 we will show that $I^{r}\cap I_{Y}$ is a cancellative
semigroup, for some $r\geq 1$.

Before proving this we  introduce some notation and develop some
machinery.  For $a,b\in S$ we write $a\tau b$ if there exists
$z\in I=S_{i-1}'\cap S_{i-1}$ so that $az=bz\not\in J$. Notice $a$
and $b$ have the same length. Hence, for a given $a\in S$, there
are only finitely many $b$ so that $a\tau b$.

Let $\cal A$ be the set of all elements $d\in I$ such that every
proper initial segment of $d$ is not in $I$. In other words, $\cal
A$ is the (unique) minimal set of generators of $I$ as a right
ideal of $S$. By the ascending chain condition on right ideals in
$S$ this is a finite set.

Let $a \tau b$ for some $a\in I, b\in S$; so $az=bz\not \in J$ for
some $z\in I$. Let $Y,Z\subseteq X$ be such that $a\in D_{Y}'$ and
$b\in D_{Z}'$.  Because $S$ is left non-degenerate,  $az\not \in
J$ implies that $z\in D_{Y}$. Since $S$ is also right
non-degenerate and because $bz\not \in J$ we thus obtain
$Z\subseteq Y$. In particular $a,b\in Sx$ for every $x\in Z$.

Choose $s\in S$ such that $a=a's, b=b's$ for some $a'\in I, b'\in
S$ such that if $a'=a''x, b=b''x$ for $x\in X, a'',b''\in S$ then
$a''\not \in I$. Notice that $a'\tau b'$. The previous paragraph
implies that $a'\in {\cal A}$.

Let $a_{j},b_{j}, j=1,\ldots ,q$, be all pairs such that
$a_{j}\tau b_{j}$, $a_{j}\neq b_{j}$ and $a_{j}\in {\cal A}$. As
remarked earlier, there are only finitely many such pairs of
elements.
 Let
$z_{j}\in I$ so that $a_{j}z_{j}=b_{j}z_{j}\not\in J$. By the
above, for every $a\in I, b\in S$ such that $a\tau b$ and $a\neq
b$ we have $a=a_{j}s, b=b_{j}s$ for some $s\in S$. Consider $Y$
that satisfies the conditions in Case~2 and assume that $a,b\in
D_{Y}'$. We claim that
\begin{eqnarray} \label{cancel}
at=bt \not \in J \;\; \mbox{ for every } t\in I^{N}\cap I_{Y}
\end{eqnarray}
where $N$ is the maximum of all $|z_{j}|$, $j=1, \ldots , q$.

So $a=a_{j}s$ and $b=b_{j}s$ for some $j$ and some $s\in S$. Since
$a_{j}\in I\setminus J$, there exist $W,Z \subseteq X$, each of
cardinality $i-1$, so that $a_{j} \in I_{WZ}$. As $a\in
D_{Y}'\setminus J$ we thus get that $a\in I_{WY}$. Moreover,
$a_{j}z_{j}=b_{j}z_{j}\not\in J$ yields that $z_{j}\in I_{ZV}$ for
some $V\subseteq X$. Now $a_{j}s=a\in D_{Y}'$ implies
$aI_{Y}=a_{j}sI_{Y}$ and, because $Y$ satisfies Case~2, the former
does not intersect $J$. In particular $sI_{Y}\subseteq I_{ZY}$.

Let $t\in I^{N}\cap I_{Y}$. Then $st\in s(I^{N}\cap
I_{Y})\subseteq D_{Z}\cap I^{|z_{j}|}\subseteq z_{j}S$ by
Lemma~\ref{claim}. So $st=z_{j}u$ for some $u\in S$. Now
$at=a_{j}st=a_{j}z_{j}u$ and $at\in I_{WY}I_{Y}$, and thus $at\not
\in J$. Similarly $bt=b_{j}st=b_{j}z_{j}u$ whence $at=bt$. This
proves the claim.

For every $Y$ that satisfies Case~2 choose $c_{Y}\in I^{N}\cap
I_{Y}$. Write $r=\max \{ |c_{Y} | +1 \}$ and let $T=I^{r}$. Then
$T/(J\cap T)=I^{r}/(J\cap I^{r})$ has a matrix pattern $T/(J\cap
T)=\bigcup _{Y,W}T_{YW}\cup \{ 0 \}$ where $T_{YW}=(I_{YW}\cap T)
\setminus J$ and $Y,W$ run through a subset of the set of
$i-1$-element subsets of $X$. The `diagonal components' are
$T_{Y}=T\cap I_{Y}$. We know that if there exist $a\in T_{YZ}$ and
$b\in T_{ZW}$ (so $Z$ satisfies Case~2) then $ab\in T_{YW}$. In
particular, if $T_{YZ}$ and $T_{ZW}$ are nonempty, then also
$T_{YW}$ is nonempty.

Let $A$ be a maximal subsemigroup of $T/(J\cap T)$ of the form
$A=\bigcup _{Y,W\in {\cal P}} T_{YW}$ $\cup \{ 0\}$ where $\cal P$
is a set of $i-1$-element subsets of $X$ such that every $T_{YW}$
is not empty. Let $Y\in {\cal P}$. Suppose that $T_{YW} \neq
\emptyset$ for some $W \not\in {\cal P}$ of cardinality $i-1$.
Then $\emptyset \neq T_{ZY}T_{YW}\subseteq T_{ZW}$ for every $Z\in
{\cal P}$. Clearly $T_{W}\neq \emptyset $ because $W$ satisfies
Case~2. Using the maximality of $\cal P$ it is now easy to see
that $B=\bigcup_{Y\in {\cal P}, \; V\not\in {\cal P}} T_{YV} \cup
\{ 0 \}$ is a right ideal of $T/(J\cap T)$. However, if $b\in B$
and $0\neq bs \in A$ for some $s\in S/J$, then $0\neq bsx\in A$
for some $x\in A$. Since $sx \in T/(J\cap T)$, it follows that
$bsx \in B$,  a contradiction. This shows that $B$ is a right
ideal of $S/J$. From the matrix pattern it follows that it is
nilpotent and this contradicts with the definition of $J$.
Consequently, $T_{YV}=\emptyset$, and similarly $T_{VY}=\emptyset$
for every $V\not \in {\cal P}$ of cardinality $i-1$. Therefore we
get a decomposition $T/(J\cap T)=A_{1}\cup \cdots \cup A_{k}$ for
some $k$, where each $A_{i}$ is of the `square type', as $A$
above. This union is $0$-disjoint, $A_{i}$ are ideals of $S/J$ and
$A_{i}A_{j}=0$ for $\neq j$. Fix some $A=A_{i}$, say $i=1$.

Let $Y$ be such that $T_{Y}\subseteq A$. Let $c=c_{Y}\in I^{N}\cap
I_{Y}$. We now prove that if $a,b\in T\cap D_{Y}'$ satisfy
$az=bz\not\in J$ for some $z\in T$ then $a=b$. As $z\in
D_{Y}\setminus J$, we must have $z\in A$. Then $z(T\cap D_{Y}')
\cap T_{Y}\neq \emptyset $, so we may assume that $z\in T_{Y}$.
Since $r\geq |c|+1$, by the dual of Lemma~\ref{claim} we may write
$a=a'c, b=b'c$ for some $a',b'\in I$. Moreover $a',b'\in D_{Y}'$
because $c\in I_{Y}$ and $a',b'\in I$ (use the generalised matrix
pattern). Then $a'cz=b'cz$ and $a'\tau b'$.  As $a'c=b'c$ by
(\ref{cancel}), we obtain $a=b$. Repeating this for every $Y$ with
$T_{Y}\subseteq A$ we show that $A$ has the property that
$az=bz\neq 0$ implies $a=b$. By a symmetric argument we may also
obtain that $A$ has the property
\begin{eqnarray}
   \mbox{if } a,b,z\in A,\; \mbox{ and }  az=bz\neq 0
 \mbox{ or } za=zb \neq 0 &\mbox { then }& a=b \label{can}
\end{eqnarray}
In particular, if $Y$ satisfies Case~2, then the diagonal
components $T_{Y}=I_{Y}\cap T$ are cancellative semigroups.

Let $Q=\{ a_{1}+\cdots +a_{m} : a_{i}\in T_{Y_{i}} \}$ where
${\cal P} = \{Y_{1},\ldots ,Y_{m}\}$. Let $Z=(T/(J\cap T)
)/(A_{2}\cup \cdots \cup A_{k})$. Then $Z$ may be identified with
$A$. Because of $(\ref{can})$ $Q$ consists of regular elements in
the algebra $K_{0}[A]$. Furthermore, the diagonal components
$T_{Y_{i}}$ form cancellative right and left Ore semigroups.
Indeed, from Lemma~\ref{claim} it follows that every two right
ideals of each $I_{Y_{i}}$ intersect nontrivially. This implies
easily that the same holds for the semigroup $T_{Y_{i}}$, and a
symmetric argument works for left ideals. It is then readily
verified that $Q$ is an Ore subset of the algebra $K_{0}[A]$. The
localization of $A$ with respect to $Q$ is an inverse semigroup
(it has a matrix pattern and each diagonal component is a group,
namely the group of quotients of the corresponding  $T_{Y}$).
Therefore $A$, and thus each $A_{i}$
 is a semigroup of generalised matrix type.
Hence $T/(J\cap T)$ has an ideal chain whose factors are of
generalised matrix type and which is determined by certain ideals
of $S$. Consider the ideal chain
 $$S_{i} \subseteq S_{i}\cup (J\cap T) \subseteq S_{i}\cup T
 \subseteq S_{i-1}\subseteq S .$$
We know that $S_{i}\cup (J\cap T)$ is nilpotent modulo $S_{i}$ and
$S_{i-1}$ is nilpotent modulo $S_{i}\cup T$. The factor
$(S_{i}\cup T)/(S_{i}\cup (J\cap T))$ is naturally identified with
$T/(J\cap T)$ because $S_{i}\cap T \subseteq J$. It follows that
$S/S_{i}$ has an ideal chain of the type described in
Theorem~\ref{sufficient}. This completes the inductive step, and
thus we have shown that $K[S]$ is right and left noetherian.

Finally, from Theorem~\ref{growth} it now follows that $K[S]$
satisfies a polynomial identity. $\Box$

\vspace{10pt} In the last paragraph of the proof we have shown
that each $A_{i}$ is an order in a completely 0-simple inverse
semigroup, in the sense of Fountain and Petrich. While this is an
easy consequence of the properties of $T$ proved before and of the
main results of \cite{fountain}, we used a simple localization
technique at the semigroup algebra level, rather than referring to
these nontrivial semigroup theoretical results.

The following is a direct consequence of the proof of
Theorem~\ref{noether}.

\begin{corollary}           \label{ideal}
Assume that $S$ is a right and left non-degenerate semigroup of
skew type. Then $S$ has a cancellative ideal $I$. Namely,
$S_{X}^{N}$ is such an ideal for some $N\geq 1$.
\end{corollary}

\section{Cancellative congruence and the prime radical}

Let $\rho$ be the least cancellative congruence on a semigroup of
skew type $S$. So it is the intersection of all congruences $\sim
$ on $S$ such that $S/\! \sim $ is cancellative. Let $\rho _{1}$
be the smallest congruence on $S$ containing all $(s,t)$ such that
$su=tu$ or $us=ut$ for some $u\in S$.  Suppose we have already
constructed $\rho _{n}$. Let $\rho _{n+1}$ be the smallest
congruence on $S$ that contains all $(s,t)$ with $(su,tu)\in \rho
_{n}$ or $(us,ut)\in \rho _{n}$ for some $u\in S$. We claim that
$\rho =\bigcup _{n\geq 1}\rho _{n}$. Indeed,  if $(su,tu)\in
\bigcup _{n\geq 1}\rho _{n} $, then $(su,tu)\in \rho _{n}$ for
some $n\geq 1$. Hence $(s,t)\in \rho _{n+1}$. It follows that
$\bigcup _{n\geq 1}\rho _{n}$ is right cancellative. Similarly, it
is left cancellative, so that $\rho \subseteq \bigcup _{n\geq
1}\rho _{n}$. For the converse first note that $\rho _{1}\subseteq
\rho$. Then, by induction one shows easily that $\rho_{n}\subseteq
\rho$ for every $n\geq 1$. Hence $\rho = \bigcup_{n\geq 1}
\rho_{n}$, as claimed.

It is easy to see (by induction) that every $\rho _{n}$ is
homogeneous, because the defining relations of $S$ are
homogeneous. It follows that $\rho $ is homogeneous.

 From now on we assume that $S$ is left and right non-degenerate.
So, by Lemma~\ref{nondegener}, $S$ satisfies: $xS\cap yS\neq
\emptyset $ for every $x,y\in S$. Define a relation $\sim $ on $S$
by: $a\sim b$ if $ax=bx$ for some $x\in S$. We claim that $\sim $
is a congruence on $S$. Suppose $a\sim b$ and $b\sim c$. Then
$ax=bx, by=cy$ for some $x,y \in S$. There exist $u,w\in S$ such
that $xu=yw$. Thus $$axu=bxu=byw=cyw=cxu$$ and so $a\sim c$. Next,
if $z\in S$, then $zs=xt$ for some $s,t\in S$. Then
$$azs=axt=bxt=bzs$$ and $az\sim bz$. It follows that $\sim $ is a
congruence on $S$. It is clear that it is the least congruence on
$S$ such that $S/ \! \sim $ is right cancellative.

\begin{lemma}          \label{ideal2}
Let $T$ be a semigroup with a cancellative ideal $J$. Assume that
$J$ has a group of quotients $G$. Define $\widehat{T}= (T\setminus
J) \cup G$. Then $\widehat{T}$ has a semigroup structure extending
that of $T$.
\end{lemma}
\noindent {\bf Proof.} The multiplication on $\widehat{T}$ is
defined by $$t(ab^{-1})=(ta)b^{-1}$$ for $a,b\in J$ and $t\in T$.
Similarly, one defines the left multiplication by elements of $G$.
Associativity can be easily checked. $\Box$

\vspace{10pt} By Lemma~\ref{ideal} there exists $N\geq 1$ such
that $I=S_{X}^{N}$ is a cancellative ideal of $S$. We know that
$I$ has a group of quotients and thus by Lemma~\ref{ideal2} we
have the semigroup $\widehat{S}=(S\setminus I) \cup II^{-1}$. Let
$e=e^{2}\in II^{-1}$. For any $a,b,x\in S$ we get that $(a-b)x=0$
implies $(a-b)xI=0$, and thus $(a-b)e=0$. Since $e$ is a central
idempotent this yields $e(a-b)=0$ and therefore $x(a-b)=0$. So, by
symmetry, we obtain  that the following conditions are equivalent:
(1) $(a-b)x=0$, (2) $x(a-b)=0$, (3) $Ix=0$ and (4) $xI=0$. It
follows that the least right cancellative congruence coincides
with the least left cancellative congruence on $S$. Note that
$\rho$ is finitely generated, as a right congruence, since $K[S]$
is right noetherian. So, we have proved the following result.

\begin{proposition}            \label{congruence}
Let $S$ be a left and right non-degenerate semigroup of skew type.
Then the least right cancellative congruence on $S$ coincides with
the least cancellative congruence on $S$ and it is defined by
$a\rho b$ if $ax=bx$ for some $x\in S$. Moreover the ideal of
$K[S]$ determined by $\rho $ is of the form $I(\rho)= \sum
_{i=1}^{k} (a_{i}-b_{i})K[S]$ for some $k\geq 1$ and
$a_{i},b_{i}\in S$.
\end{proposition}

Recall that by definition $I(\rho )$ is the kernel of the natural
homomorphism $K[S]\longrightarrow K[S/\rho ]$. The congruence
$\rho$ is actually important for the description of the prime
radical $B(K[S])$ of $K[S]$. As in the proof of
Theorem~\ref{growth}, we get that $S/\rho $ has an
abelian-by-finite group of quotients. Moreover, if $\ch (K)=0$,
then $K[S/\rho ]$ is semiprime (see for example \cite{book},
Theorem~7.19). In particular $B(K[S])\subseteq I(\rho )$, the
ideal of $K[S]$ determined by $\rho $.

We have seen that $a\rho b$ if and only if $(a-b)I=0=I(a-b)$. So
$I(\rho )I=II(\rho )=0$. Hence, if $P$ is a prime ideal of $K[S]$
with $P\cap S = \emptyset$ then $I(\rho )\subseteq P$. If, on the
other hand, $P$ is a prime ideal with $P\cap S\neq \emptyset$,
then there exists $b\in P\cap S_{X}$. So, by Lemma~\ref{claim},
$S_{X}^{k}=bS \subseteq P$ for some positive integer $k$. It
follows that $S_{X}\subseteq P$.

Suppose that $\alpha \in K[S]$ belongs to the left annihilator
$\ann_{l}(I)$ of $I$ in $K[S]$. Then $\alpha s =0$ for all $s\in
I$. Write $\alpha =\alpha_{1} +\cdots + \alpha_{m}$ with $|\supp
(\alpha_{i})s|=1$ and $\supp (\alpha_{i})s\neq \supp
(\alpha_{j})s$ for $i\neq j$. It follows that $\alpha_{i} s=0$ for
all $i$. So the augmentation of $\alpha _{i}$ is zero and it is
clear that $\alpha _{i}\in I(\rho )$. From all the above it
follows that $I(\rho )\subseteq \ann_{l}(I) \subseteq I(\rho )$.
Hence $I(\rho )=\ann _{l}(I)$. By symmetry, $I(\rho )= \ann (I)$,
the two-sided annihilator of $I$. If $\ch (K)=0$, then $I(\rho )$
is a semiprime ideal and thus $I(\rho ) =\bigcap_{P,
I(\rho)\subseteq P} P$. So we have proved the following result. By
$X^{0}(K[S])$ we denote  the set of all the minimal primes of
$K[S]$.

\begin{proposition}          \label{radical}
If $S$ is a left and right non-degenerate semigroup of skew type,
then
\begin{enumerate}
\item
$I(\rho )=\ann (S_{X}^{N})$ for some $N\geq 1$,
\item   $I(\rho) \subseteq P$ for any $P\in X^{0}(K[S])$ with
 $P\cap S =\emptyset$,
\item $S_{X}\subseteq P$ for any $P\in X^{0}(K[S])$ with
$P\cap S \neq \emptyset$.
\end{enumerate}

If, furthermore, $\ch (K)=0$ then $$B(K[S])=I(\rho) \cap
\bigcap_{P\in X^{0}(K[S]), P\cap S \neq \emptyset } P = I(\rho )
\cap (\bigcap_{P\in X^{0}(K[S]),S_{X}\subseteq P} P) .$$
\end{proposition}

Note that if $S$ is a left and right non-degenerate semigroup of
skew type then there is at least one minimal prime $P$ so that
$P\cap S = \emptyset$. Indeed for otherwise the proposition
implies that $S_{X}\subseteq B(K([S])$. This yields a
contradiction as $S_{X}$ is not nil.

\section{Examples}

Our first example shows that $K[S]$ can be a noetherian PI algebra
even if $S$ is not non-degenerate.

\begin{example} Let $S=\langle
x_{1},x_{2},x_{3},x_{4}\rangle $
 be the semigroup of skew type defined by
$$x_{3}x_{2}=x_{1}x_{4} \, \mbox{ and } \, x_{4}x_{1}=x_{2}x_{3}$$
with all the remaining relations in $S$ of the form $xy=yx$. Then,
for every field $K$, $K[S]$ is a noetherian PI-algebra. Moreover
$B(K[S])=I(\rho )$, and $I(\rho )$ is the only minimal prime of
$K[S]$.
\end{example}

\noindent {\bf Proof.} By the defining relations we get
$$x_{3}x_{3}x_{2}=x_{3}x_{1}x_{4}=x_{1}x_{3}x_{4}=x_{1}x_{4}x_{3}=
x_{3}x_{2}x_{3}=$$
$$=x_{3}x_{4}x_{1}=x_{4}x_{3}x_{1}=x_{4}x_{1}x_{3}=
x_{2}x_{3}x_{3}.$$ So $x_{3}^{2}\in Z(S)$.  Similarly
$x_{i}^{2}\in Z(S)$ for $i=1,2,3,4$. It is easy to see that $S=\{
x_{1}^{a_{i}}x_{2}^{a_{2}}x_{3}^{a_{3}}x_{4}^{a_{4}} | a_{i}\geq
0\}$. It follows that $K[S]$ is a finite module over $K[A]$, where
$A=\langle x_{1}^{2},x_{2}^{2},x_{3}^{2},x_{4}^{2} \rangle $.
Therefore $K[S]$ is right and left noetherian. Notice that $S$ is
right and left degenerate.

We claim that $B(K[S])=I(\rho )$ and the least cancellative
congruence $\rho $ on $S$ coincides with the congruence determined
by the natural homomorphism $\phi : K[S]\longrightarrow K[C]$,
where $C$ is the commutative monoid obtained from $S$ by adding
all the commutator relations to the defining relations of $S$. So
$C=\langle a_{1},a_{2},a_{3},a_{4} \mid a_{i}a_{j}=a_{j}a_{i},
a_{1}a_{4}=a_{3}a_{2} \rangle $. In fact, we have seen above that
$$x_{3}(x_{2}x_{3}-x_{3}x_{2})=0=(x_{2}x_{3}-x_{3}x_{2})x_{3}.$$
Similarly one shows that
$$y(x_{2}x_{3}-x_{3}x_{2})=0=(x_{2}x_{3}-x_{3}x_{2})y$$ for every
$y\in \{ x_{1},x_{2},x_{3},x_{4}\}$. Also,
$$y(x_{1}x_{4}-x_{4}x_{1})=0=(x_{1}x_{4}-x_{4}x_{1})y.$$ So $\ker
(\phi ) \subseteq I(\rho )$. Since $C$ embeds in a torsion-free
group, $K[C]$ is a domain and we get $I(\rho )= \ker (\phi )$.
Also $B(K[S])\subseteq \ker (\phi)$, while $(\ker \phi )^{2}=0$ by
the displayed formulas. It follows that $B(K[S])=I(\rho )$. In
particular, $I(\rho )$ is the only minimal prime of $K[S]$, so
$K[S]$ has no minimal primes intersecting $S$. $\Box$

\vspace{10pt}

The following example shows that a right non-degenerate semigroup
$S$ of skew type  does not always yield a right noetherian algebra
$K[S]$.

\begin{example}
Let $S=\langle x_{1},x_{2},x_{3} \rangle $ be the monoid defined
by the relations $$ x_{2}x_{1}=x_{3}x_{1}, x_{1}x_{2}=x_{3}x_{2},
        x_{1}x_{3}=x_{2}x_{3} .$$
Then $S$ is right non-degenerate but not left non-degenerate and
$K[S]$ is neither right nor left noetherian. Furthermore, $S$ is
left cancellative and $\GK (K[S])=2$.
\end{example}
\noindent {\bf Proof.} Consider the elements
$a_{n}=x_{1}x_{2}^{n+1}-x_{1}^{n}x_{2}^{2}\in K[S]$, $n=2,3 \ldots
$. We claim that $a_{n}K[S]= \lin_{K} \{ a_{n}x_{2}^{j} | j\geq 0
\}$. Indeed, first note that
$$x_{1}x_{2}x_{1}=x_{1}x_{3}x_{1}=x_{2}x_{3}x_{1}=x_{2}x_{2}x_{1}.$$
Next
\begin{eqnarray*} x_{1}x_{2}x_{2}x_{1}=x_{1}x_{2}x_{3}x_{1}=
x_{1}x_{1}x_{3}x_{1}=x_{1}x_{1}x_{2}x_{1}=x_{1}x_{3}x_{2}x_{1}= \\
x_{2}x_{3}x_{2}x_{1}=x_{2}x_{1}x_{2}x_{1}=x_{2}x_{1}x_{3}x_{1}=
x_{2}x_{2}x_{3}x_{1}=x_{2}x_{2}x_{2}x_{1}. \end{eqnarray*} Using
induction, we then also get for every $a>2$
$$x_{1}x_{2}^{a}x_{1}=x_{1}x_{2}^{a-1}x_{3}x_{1}=
x_{1}x_{2}^{a-2}x_{1}x_{3}x_{1}=x_{2}^{a-1}x_{1}x_{3}x_{1}=
x_{2}^{a}x_{3}x_{1}=x_{2}^{a+1}x_{1} .$$ Now, for every $n\geq 2$
we get $$a_{n}x_{1}=0$$ and $$a_{n}x_{3}=
(x_{1}x_{2}^{n}-x_{1}^{n}x_{2})x_{2}x_{3}=
(x_{1}x_{2}^{n}-x_{1}^{n}x_{2})x_{1}x_{3}=0 .$$ So the claim
follows.

Now each element $x_{1}^{k}x_{2}^{q}$ can only be rewritten as
$syx_{2}^{q}$ for some $s\in S$ and $y\in \{ x_{1},x_{3} \}$. It
then easily follows that for $n\geq 3$ there do not exist
$\lambda_{j}\in K$ so that
  $$a_{n}=\sum _{j=2}^{n-1} \lambda_{j}
    (x_{1}x_{2}^{j+1}-x_{1}^{j}x_{2}^{2})x_{2}^{n-j} .$$
Therefore, $a_{n}\not \in \sum _{j=2}^{n-1}a_{j}K[S]$ for every
$n$. So, indeed $K[S]$ is not right noetherian (however, $S$
satisfies the ascending chain condition by Proposition~\ref{acc}).

If $k<n$ then $x_{2}x_{1}^{n} \not \in Sx_{2}x_{1}^{k}$. This is
clear from the defining relations. Namely, for every $s\in S$ the
element $sx_{2}x_{1}^{k}$ can only be rewritten in the form
$tx_{2}x_{1}^{k}$ or $tx_{3}x_{1}^{k}$ for some $t\in S$. So $S$
does not satisfy the ascending chain condition on left ideals and
$K[S]$ is not left noetherian.

It can be verified that $\GK (K[S])=2$. From the relations it also
follows easily that $S$ is left cancellative. $\Box$

\vspace{10pt} Our third example satisfies the cyclic condition,
but the defining relations do not yield a Gr\"{o}bner basis, so it
is not of binomial type studied in \cite{gateva},\cite{binom}. The
aim is to show that one can get important structural information
on $K[S]$. In particular we determine all minimal primes and the
prime radical of $K[S]$. Recall that $K[S]$ is an affine PI
algebra which is left and right noetherian by Theorem~\ref{growth}
and Proposition~\ref{cyclic}.

\begin{example} \label{exA1} Let $S=\langle x_{1},x_{2},x_{3},x_{4}\rangle $
be given by the presentation
 $$
  x_{4}x_{3}=x_{1}x_{4} ,  x_{4}x_{2}=x_{2}x_{4} ,
   x_{4}x_{1}=x_{3}x_{4},$$
 $$   x_{3}x_{2}=x_{1}x_{3} ,
   x_{3}x_{1}=x_{2}x_{3} ,   x_{2}x_{1}=x_{1}x_{2} .$$
The minimal primes of $K[S]$ are the ideals $P_{1}=(x_{1}-x_{2},
x_{2}-x_{3})=I(\rho )$, $P_{2}=(x_{4})$, $P_{3}=(x_{1},x_{3})$ and
$P_{4}=(x_{2})$. Moreover, $K[S]$ is semiprime, has dimension
three and $$S=\langle x_{1},x_{2},x_{3} \rangle \cup \langle
x_{1}\rangle \langle x_{4}\rangle \cup \langle x_{2}\rangle
\langle x_{4} \rangle \cup \langle x_{3} \rangle \langle x_{4}
\rangle \cup \langle x_{1}\rangle x_{2} \langle x_{4} \rangle .$$
\end{example}

\noindent {\bf Proof.} First note that the following equalities
hold in $S$:
 $$
 x_{1}x_{3}x_{4}=x_{1}x_{4}x_{1}=x_{4}x_{3}x_{1}=x_{4}x_{2}x_{3}
   = x_{2}x_{4}x_{3}=x_{2}x_{1}x_{4}=x_{1}x_{2}x_{4}
 \;\;\; \; (7)$$
and
 $$
  x_{1}x_{2}x_{4}=x_{1}x_{4}x_{2}=x_{4}x_{3}x_{2}
  =
  x_{4}x_{1}x_{3}=x_{3}x_{4}x_{3}=x_{3}x_{1}x_{4}=x_{2}x_{3}x_{4}.
  \; \;\;\;(8)$$
So  $x_{1}x_{3}x_{4}=x_{1}x_{2}x_{4}$ and $x_{2}x_{1}x_{4}=
x_{2}x_{3}x_{4}$. Therefore $P_{1}=(x_{2}-x_{3},\; x_{1}-x_{3})
\subseteq I(\rho )$. As $K[S]/P_{1} \cong K[Y_{1},Y_{4}]$, a
polynomial ring in two commuting variables, we get that $P_{1}$ is
a prime ideal of $K[S]$. So, by Proposition~\ref{radical} and its
following remark, $P_{1}$ is a minimal prime ideal of $K[S]$ (it
has depth $2$), $I(\rho )=P_{1}$,  and $P_{1}$ is the only minimal
prime of $K[S]$ intersecting $S$ trivially.

Second note that $x_{4}$ is a normalizing element of $S$ and thus
also a normalizing element of $K[S]$. Also $K[S]/(x_{4})\cong
K[\langle x_{1},x_{2},x_{3} \rangle ]$ and because $\langle
x_{1},x_{2},x_{3}\rangle$ is a binomial semigroup, we get that
$(x_{4})$ is a prime ideal of depth $3$.

 Now, suppose $P$ is a prime ideal
of $K[S]$ that does not contain $x_{4}$. The equations (7) and (8)
yield that
 $$I=(x_{1}(x_{2}-x_{3}),(x_{1}-x_{2})x_{3}) \subseteq P.$$
In the classical ring of quotients $Q_{cl}(K[S]/P)$ the element
$\overline{x_{4}}$ is invertible  (as it is regular in $K[S]/P$)
and this element acts via conjugation on the set $\{
\overline{x_{1}},\overline{x_{2}}, \overline{x_{3}}\}$. Applying
this conjugation action on the equations
$\overline{x_{1}}\overline{x_{2}}=\overline{x_{1}}\overline{x_{3}}=
\overline{x_{2}}\overline{x_{3}}$ yields
$\overline{x_{3}}\overline{x_{2}}=\overline{x_{3}}\overline{x_{1}}=
\overline{x_{2}}\overline{x_{1}}$. As $x_{1}x_{2}=x_{2}x_{1}$ and
thus
$\overline{x_{1}}\overline{x_{2}}=\overline{x_{2}}\overline{x_{1}}$
we get that the monoid $\langle
\overline{x_{1}},\overline{x_{2}},\overline{x_{3}}\rangle$ is
abelian.  It is easily verified that $\langle
\overline{x_{1}},\overline{x_{2}},\overline{x_{3}}\rangle =
\langle \overline{x_{1}} \rangle \cup  \langle \overline{x_{2}}
\rangle \cup \langle \overline{x_{3}} \rangle \cup \langle
\overline{x_{1}}\rangle \overline{x_{2}}$. It follows that
$K[S]/P$ is an epimorphic image of $K[S/\tau ]$, where $\tau$ is
the smallest congruence generated by the relations in $S$ and the
extra relations $x_{1}x_{2}=x_{1}x_{3}=x_{3}x_{1}=
x_{2}x_{3}=x_{3}x_{2}$. Denote the image of $x_{i}$ in $S/\tau$ by
$y_{i}$. Then we get $S/\tau =\left( \langle y_{1} \rangle \cup
\langle y_{2} \rangle \cup \langle y_{3} \rangle \cup \langle
y_{1}\rangle y_{2}\right) \langle y_{4} \rangle $. Moreover
$y_{4}$ acts on $T=\langle y_{1},y_{2},y_{3}\rangle \subseteq
S/\tau$ via an automorphism $\sigma$ of finite order. It follows
that $K[S/\tau ] =(K[T])[y_{4},\sigma ]$, a skew polynomial ring.
Now the commutative semigroup $T$ is a semilattice of cancellative
semigroups, each yielding an algebra which is a domain. Hence
$K[T]$ is semiprime and thus so is the skew polynomial ring
$K[S/\tau ]$. Moreover,
 $$y_{1}(y_{2}-y_{3})=(y_{1}-y_{2})y_{3}=0.$$
Thus if $Q$ is a minimal prime ideal of the abelian algebra $K[T]$
then $Q$ contains one of the following ideals:
 $$(y_{1},y_{3}),\; (y_{1},y_{1}-y_{2})=(y_{1},y_{2}),\;
 (y_{3},y_{2}-y_{3})=(y_{3},y_{2}), \;
 (y_{2}-y_{3},y_{1}-y_{3}).$$
It is easily seen that each of these ideals is a prime ideal of
$K[T]$ of depth $1$.  Hence these are all the minimal
 prime ideals  of  $K[T]$. Under the action
of $\sigma$ there are thus precisely three orbits of minimal
primes in $K[T]$. Hence the minimal $\sigma$-primes of $K[T]$ are
$(y_{1},y_{3})$, $(y_{1},y_{2})\cap (y_{2},y_{3})=(y_{2})$ and
$(y_{1}-y_{3},y_{2}-y_{3})$. Note also that
 $$\hspace*{2cm}
   T\cap
  (y_{4})\cap (y_{2})\cap (y_{1},y_{3})= \{ y_{1}^{\alpha} y_{2}
   y_{4}^{\gamma}\mid \alpha , \gamma >0\}. \hspace*{2.5cm} (9)
  $$

It is easily seen (and well known from standard results on ${\bf
Z}$-graded rings) that the minimal primes of the skew polynomial
algebra $K[S/\tau ]=(K[T])[y_{4},\sigma]$ are all ideals of the
type $M[y_{4},\sigma ]$ with $M$ a minimal $\sigma$-prime ideal of
$K[T]$. Therefore the minimal primes of $K[S/\tau ]$ are
 $(y_{2})$, $(y_{1},y_{3})$ and $(y_{1}-y_{3},y_{2}-y_{3})$.

All the above implies that if  $P$ is  a
 prime ideal of $K[S]$ that does not
contain $x_{4}$, then $P$ contains one of the following
incomparable prime ideals of depth $2$:
 $$J+ (x_{2})= (x_{2}) =P_{4}, \;
 J+ (x_{1},x_{3})=(x_{1},x_{3})=P_{3} \mbox{ or } P_{1},$$
where $J$ is the kernel of the natural epimorphism
$K[S]\longrightarrow K[S/\tau ]$. As all these primes are
incomparable with the prime $P_{2}=(x_{4})$, we get that indeed
$P_{1},P_{2},P_{3},P_{4}$ are all the minimal prime ideals of
$K[S]$. Because $P_{2}$ has maximal depth, $K[S]$ has dimension
$3$.

 From (9) we get that
 $$(x_{4}) \cap (x_{2})\cap (x_{1},x_{3}) \subseteq J+K[\{ x_{1}^{\alpha} x_{2}
   x_{4}^{\gamma} \mid \alpha , \gamma >0\} ]\subseteq
   (x_{4}) \cap (x_{2}) \cap (x_{1},x_{3}).$$
Hence it follows easily that
 $S\cap (x_{4})\cap (x_{2}) \cap (x_{1},x_{3})= \{ x_{1}^{\alpha} x_{2}
   x_{4}^{\gamma} \mid \alpha , \gamma >0\} .$

Since $K[S]/(x_{1}-x_{2},x_{2}-x_{3})\cong K[Y_{1},Y_{4}]$, a
polynomial algebra in commuting variables, it also easily follows
that $x_{1}^{\alpha}x_{2}x_{4}^{\gamma}=x_{1}^{\alpha
'}x_{2}x_{4}^{\gamma'}$ if and only if $\alpha =\alpha '$ and
$\gamma =\gamma '$. Hence $(x_{1}-x_{2},x_{2}-x_{3}) \cap K[\{
x_{1}^{\alpha}x_{2}x_{4}^{\gamma}\mid \alpha , \gamma > 0\}] =\{
0\}$. As $P_{2}\cap P_{3}\cap P_{4} = K[\{ x_{1}^{\alpha} x_{2}
   x_{4}^{\gamma} \mid \alpha , \gamma >0\}]$ we  indeed obtain
that $K[S]$ is semiprime.

Now earlier we have shown that $x_{1}x_{2}x_{4}=x_{1}x_{3}x_{4}=
x_{2}x_{3}x_{4}$. Using these and the defining relations for $S$,
it is easy to check that $x_{1}^{\alpha} x_{2}^{\beta}
x_{3}^{\gamma} x_{4}=x_{1}^{\alpha +\beta + \gamma -1} x_{2}
x_{4}$ in case that at least two of the exponents $\alpha , \beta
,\gamma$ are nonzero. Hence $$S=\langle x_{1} , x_{2},x_{3}\rangle
\cup \{ x_{1}^{\alpha}x_{4}^{\gamma},
x_{2}^{\alpha}x_{4}^{\gamma}, x_{3}^{\alpha}x_{4}^{\gamma},
x_{1}^{\alpha}x_{2}x_{4}^{\gamma}\mid \alpha , \gamma \geq 0\} .$$
$\Box$

\vspace{10pt}

\noindent \begin{tabular}{llllllll}
 T. Gateva-Ivanova && E. Jespers  \\ Institute of Mathematics
 && Department of
Mathematics \\ Bulgarian Academy of Sciences && Vrije Universiteit
Brussel \\ Sofia 1113 && Pleinlaan 2  \\
 Bulgaria && 1050 Brussel, Belgium  \\

 \\

  J. Okni\'{n}ski   \\
  Institute of Mathematics  \\
  Warsaw University\\
  Banacha 2  \\
  02-097 Warsaw, Poland
\end{tabular}
\end{document}